\documentclass[10pt]{article}
\font\smallit=cmti10

\usepackage{mathtools,amssymb,amsthm}
\usepackage{float, graphicx}
\usepackage[style=alphabetic, backend=bibtex]{biblatex}
\usepackage{hyperref}

\newcommand{\bbn}{\mathbb{N}}
\newcommand{\bbr}{\mathbb{R}}

\numberwithin{equation}{section}

\newtheorem{defn}{Definition}
\newtheorem{thm}{Theorem}

\makeatletter

\renewcommand\section{\@startsection {section}{1}{\z@}
	{-30pt \@plus -1ex \@minus -.2ex}
	{2.3ex \@plus.2ex}
	{\normalfont\normalsize\bfseries}}

\renewcommand\subsection{\@startsection{subsection}{2}{\z@}
	{-3.25ex\@plus -1ex \@minus -.2ex}
	{1.5ex \@plus .2ex}
	{\normalfont\normalsize\bfseries}}

\renewcommand{\@seccntformat}[1]{\csname the#1\endcsname. }

\makeatother

\begin{document}
	\begin{center}
		\textbf{\uppercase{Spacing Distribution of a Bernoulli Sampled Sequence}}
		\vskip 20pt
		\textbf{Abigail L. Turner}\\
		{\smallit Department of Mathematics, Cornell University, Ithaca, NY 14853}\\
		\texttt{alt86@cornell.edu}
		\vskip 10pt
		\textbf{Ananya Uppal}\\
		{\smallit Department of Mathematical Sciences, Carnegie Mellon University,  Pittsburgh, PA 15213}\\
		\texttt{auppal@andrew.cmu.edu}
		\vskip 10pt
		\textbf{Peng Xu}\\
		{\smallit Department of Industrial and Enterprise Systems Engineering, University of Illinois, Urbana, IL 61801}\\
		\texttt{pengxu1@illinois.edu}
	\end{center}
	\vskip 30pt
	
	\centerline{\bf Abstract}
	
	\noindent
	We investigate the spacing distribution of sequence \[S_n=\left\{0,\frac{1}{n},\frac{2}{n},\dots,\frac{n-1}{n},1\right\}\] after Bernoulli sampling. We describe the closed form expression of the probability mass function of the spacings, and show that the spacings converge in distribution to a geometric random variable.
	
	\pagestyle{myheadings}
	\baselineskip=12.875pt
	\vskip 30pt 
	
	\section{Introduction}
	Given \(n\in\bbn\), consider the sequence \[S_n=\left\{0,\frac{1}{n},\frac{2}{n},\dots,\frac{n-1}{n},1\right\}.\] If we sample \(S_n\) in a manner such that each element has a constant probability of survival \(p\), a new sequence \(S'_n\subseteq S_n\) will be obtained, and \(|S'_n|\) will have binomial distribution with parameters \(n+1\) and \(p\). By assuming \(S'_n\) contains more than \(i\) elements, we are able to define the random spacing:
	\begin{defn}
		The \(i\)-th \textbf{random spacing} of sequence \(S_n\) is \[\Delta_{i,n}\coloneqq s'_{i+1}-s'_i,\] where \(s'_i\in S'_n\).
	\end{defn} 
	Since we will be considering the limiting behavior of random spacing, it is also useful to define the scaled random spacing:
	\begin{defn}
		The \(i\)-th \textbf{scaled random spacing} of sequence \(S_n\) is \[D_{i,n}\coloneqq n\Delta_{i,n}.\]
	\end{defn}
	We shall study the distribution of the random spacing, especially its asymptotic behavior. More specifically, we shall provide proofs for the following theorems:
	\begin{thm}\label{thm:1}
		The probability mass function (p.m.f.) of random spacing \(\Delta_{i,n}\) is \[f_{\Delta_{i,n}}(\delta)=\frac{\displaystyle p^{i+1}(1-p)^{n\delta-1}\sum_{j=i-1}^{n-n\delta}\binom{j}{i-1}(1-p)^{j-i+1}}{\displaystyle 1-\sum_{k=0}^{i}\binom{n+1}{k}p^k(1-p)^{n+1-k}},\qquad\delta\in\left\{\frac{1}{n},\frac{2}{n},\dots,\frac{n-1}{n},1\right\}.\]
		
		In addition, the p.m.f.\ of the scaled random spacings \(D_{i,n}\) is \[f_{D_{i,n}}(d)=\frac{\displaystyle p^{i+1}(1-p)^{d-1}\sum_{j=i-1}^{n-d}\binom{j}{i-1}(1-p)^{j-i+1}}{\displaystyle 1-\sum_{k=0}^{i}\binom{n+1}{k}p^k(1-p)^{n+1-k}},\qquad d\in\left\{1,2,\dots,n-1,n\right\}.\]
	\end{thm}
	\begin{thm}\label{thm:2}
		The sequence of random variables \(\{D_{1,n}\}\), with support \(d=\{1,2,\dots\}\), converges in distribution to a geometric random variable with parameter \(p\). That is, \[\lim_{n\to\infty}F_{D_{1,n}}(d)=1-(1-p)^d.\] 
	\end{thm}
	\begin{thm}\label{thm:3}
		The sequence of random variables \(\{D_{i,n}\}\) converges in distribution to a geometric random variable with parameter \(p\) for all \(i\). That is, \[\lim_{n\to\infty}F_{D_{i,n}}(d)=1-(1-p)^d.\]
	\end{thm}
	
	\subsection{Background}
	This paper originated from \href{http://www.math.illinois.edu/igl/}{Illinois Geometry Lab} research project \emph{Random Discrete Sets} (formerly named \emph{Visibility in Random Forests}). The original goal was to study the spacing distribution of Farey fractions, whose density are given by \autocite{hall70}. We investigated the limiting behavior of said distribution, as well as the spacing distribution of randomly chosen subsets of Farey fractions, during which we discovered that the spacings of subsets of Farey fractions, obtained by Bernoulli sampling, are approximately exponentially distributed. Similar observations are later noted for other sequences on \([0,1]\), such as \(\{n\alpha\}\bmod{1}\), where \(\alpha\) is any irrational. In order to explain this phenomenon, we introduced and studied \(S_n\), a sequence with elementary structure.
	\begin{figure}[H]
		\centering
		\includegraphics[scale = 0.6]{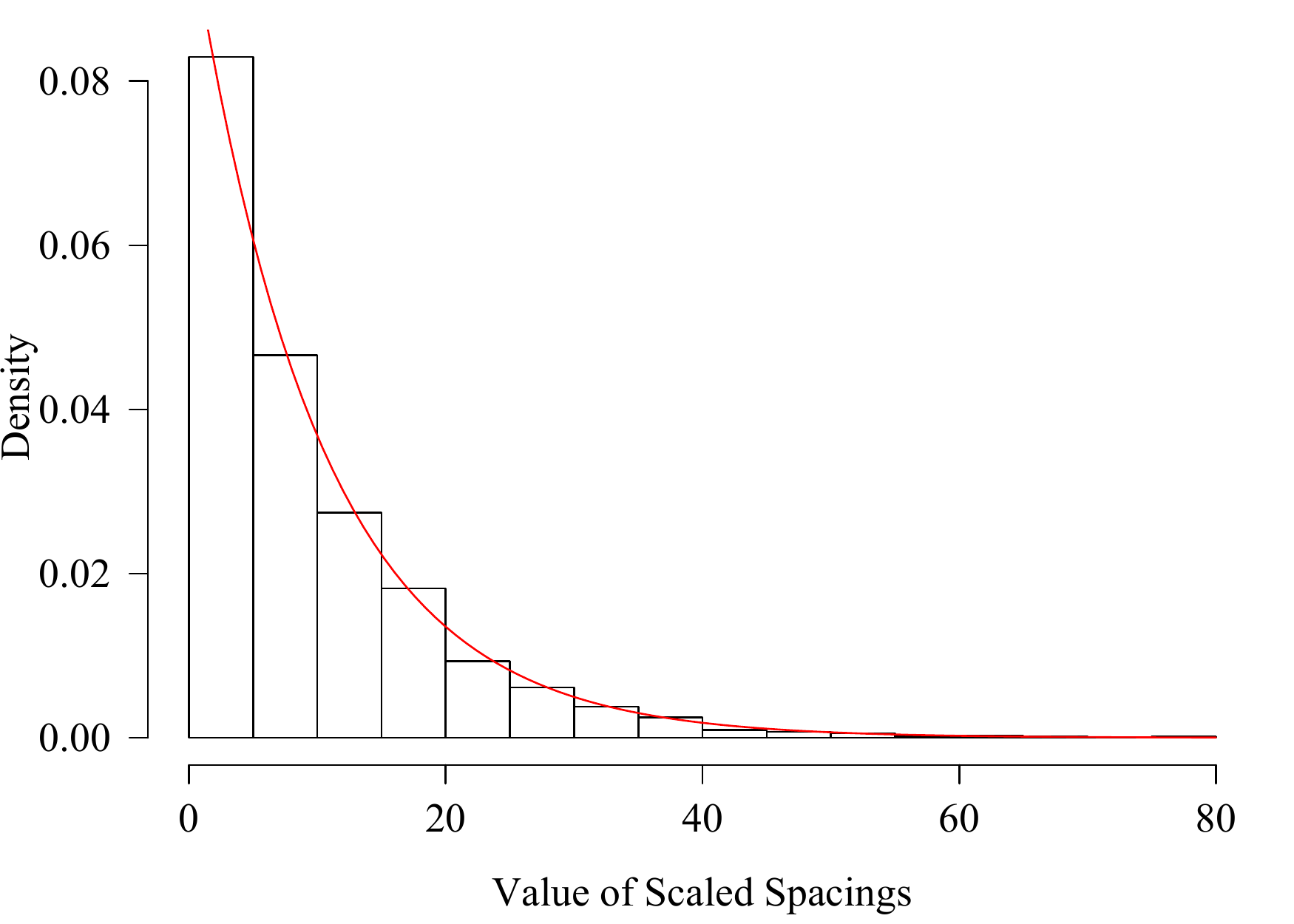}
		\caption{Density of scaled spacings of \(S_n'\), with \(n = 50000\) and \(p = 1/10\), overlapped with the density function of a exponential distribution, \(\lambda = 1/10\).}
	\end{figure}
	
	Other works that are not directly linked to this paper, but are crucial to this research include, but not limited to: \autocite{ac14}, \autocite{bhs13}, \autocite{gkcdh68}, \autocite{pyke65}, \autocite{pyke72}, and \autocite{vr88}.
	
	\subsection{Organization}
	The rest of the paper is organized as follows: in section \ref{sec:pmf}, we will derive the p.m.f.\ for the random spacings, thus proving Theorem \ref{thm:1}; in section \ref{sec:sc1}, we will derive the limiting distribution of the scaled random spacing, for the case \(n=1\), thus proving Theorem \ref{thm:2}; in section \ref{sec:gc}, we will again derive the limiting distribution of the scaled random spacings, this time for general case, thus proving Theorem \ref{thm:3}; and in section \ref{sec:ap}, we will provide an alternative proof to Theorem \ref{thm:3} from a stochastic point of view.
	
	\section{Probability Mass Function of Random Spacings}
	\label{sec:pmf}
	In this section, we will derive the p.m.f.\ for \(\Delta_{i,n}\) and \(D_{i,n}\), thereby proving Theorem \ref{thm:1}.
	
	\begin{proof}
		The p.m.f.\ of \(\Delta_{i,n}\) is defined as
		\begin{equation}\label{eq:2.1}
			f_{\Delta_{i,n}}(\delta)\coloneqq P\big(\Delta_{i,n}=\delta\mid|S'_n|>i\big),\qquad \delta\in\left\{\frac{1}{n},\frac{2}{n},\dots,\frac{n-1}{n},1\right\}.
		\end{equation}
		Using Bayes' theorem, right hand side of equation \eqref{eq:2.1} can be rewritten as the following:
		\begin{equation}\label{eq:2.2}
			P\left(\Delta_{i,n}=\delta\mid|S'_n|>i\right)=\frac{\displaystyle P\left(\Delta_{i,n}=\delta\right)\cdot P\left(|S'_n|>i\mid\Delta_{i,n}=\delta\right)}{P(|S'_n|>i)}
		\end{equation}
		To simplify equation \eqref{eq:2.2}, first rewrite \(P(\Delta_{i,n}=\delta)\) as
		\begin{equation}\label{eq:2.3}
			P\left(\Delta_{i,n}=\delta\right)=\sum_{j=0}^{n-1}P\left(s'_{i+1}-s'_i=\delta\;\middle|\;s'_i=\frac{j}{n}\right)P\left(s'_i=\frac{j}{n}\right).
		\end{equation}
		Having spacing of value \(\delta\) requires \(n\delta-1\) elements being eliminated in between \(s'_i\) and \(s'_{i+1}\), hence
		\begin{equation*}
			P\left(s'_{i+1}-s'_i=\delta\;\middle|\;s'_i=\frac{j}{n}\right)=
			\begin{dcases}
				(1-p)^{n\delta-1}p, & j\leq n-n\delta,\\
				0, & \text{else.}
			\end{dcases}
		\end{equation*}
		The number of surviving elements before the \(i\)-th surviving element has binomial distribution, thus
		\begin{equation*}
			P\left(s'_i=\frac{j}{n}\right)=
			\begin{dcases}
				\binom{j}{i-1}p^{i}(1-p)^{j-i+1}, & j\geq i-1,\\
				0, & \text{else.}
			\end{dcases}
		\end{equation*}
		Therefore, equation \eqref{eq:2.3} can be written as
		\begin{equation}\label{eq:2.4}
			P(\Delta_{i,n}=\delta)=p^{i+1}(1-p)^{n\delta-1}\sum_{j=i-1}^{n-n\delta}\binom{j}{i-1}(1-p)^{j-i+1}.
		\end{equation}
		Moreover, since \(|S'_n|\) has binomial distribution, we have
		\begin{equation}\label{eq:2.5}
			P(|S'_n|>i)=1-P(|S'_n|\leq i)=1-\sum_{k=0}^{i}\binom{n+1}{k}p^k(1-p)^{n+1-k}.
		\end{equation}
		Finally, the existence of $\Delta_{i,n}$ implies that there are at least \(i+1\) elements in \(S'_n\), hence we get
		\begin{equation}\label{eq:2.6}
			P(|S'_n|>i\mid\Delta_{i,n}=\delta)=1.
		\end{equation}
		By combining equations \eqref{eq:2.4}, \eqref{eq:2.5}, and \eqref{eq:2.6}, together with change of variable, the p.m.f.'s of both \(\Delta_{i,n}\) and \(D_{i,n}\) are obtained:
		\begin{align*}
			f_{\Delta_{i,n}}(\delta)&=\frac{\displaystyle P\left(\Delta_{i,n}=\delta\right)}{P(|S'_n|>i)}\\
			&\quad\implies f_{\Delta_{i,n}}(\delta) =\displaystyle\frac{ p^{i+1}(1-p)^{n\delta-1}\displaystyle\sum_{j=i-1}^{n-n\delta}\binom{j}{i-1}(1-p)^{j-i+1}}{1-\displaystyle\sum_{k=0}^{i}\binom{n+1}{k}p^k(1-p)^{n+1-k}}.
		\end{align*}
	\end{proof}
	
	\section{Limiting Distribution of Random Gaps}
	\numberwithin{equation}{subsection}
	In this section, we will derive the limiting distribution of \(D_{i,n}\) with respect to \(n\), thereby proving Theorem \ref{thm:2} and Theorem \ref{thm:3}.
	\subsection{Special Case: \(i=1\)}
	\label{sec:sc1}
	
	\begin{proof}
		Let \(i=1\), the p.m.f.\ of \(\Delta_{i,n}\) is 
		\[
		f_{\Delta_{i,n}}(\delta) =\displaystyle\frac{ p^{i+1}(1-p)^{n\delta-1}\displaystyle\sum_{j=i-1}^{n-n\delta}\binom{j}{i-1}(1-p)^{j-i+1}}{1-\displaystyle\sum_{k=0}^{i}\binom{n+1}{k}p^k(1-p)^{n+1-k}}
		\]
		Thus, $f_{\Delta_{1,n}}(\delta)$ is
		\[f_{\Delta_{1,n}}(\delta)=\frac{\displaystyle p^{2}(1-p)^{n\delta-1}\sum_{j=0}^{n-n\delta}(1-p)^{j}}{\displaystyle 1-(1-p)^{n+1}-(n+1)p(1-p)^n},\qquad\delta\in\left\{\frac{1}{n},\frac{2}{n},\dots,\frac{n-1}{n},1\right\}.\] 
		Because it is a geometric series,
		\[\sum_{j=0}^{n-n\delta}(1-p)^{j}=\frac{1-(1-p)^{n-n\delta+1}}{p},\] therefore, \[f_{\Delta_{1,n}}(\delta)=\frac{\displaystyle p(1-p)^{n\delta-1}(1-(1-p)^{n-n\delta+1})}{\displaystyle 1-(1-p)^{n+1}-(n+1)p(1-p)^n},\qquad\delta\in\left\{\frac{1}{n},\frac{2}{n},\dots,\frac{n-1}{n},1\right\}.\] 
		With change of variable \(D_{1,n}=n\Delta_{1,n}\) and \(d=n\delta \), an equivalent expression is 
		\[f_{D_{1,n}}(d)=\frac{\displaystyle p(1-p)^{d-1}(1-(1-p)^{n-d+1})}{1-(1-p)^{n+1}-(n+1)p(1-p)^n},\qquad d\in\left\{1,2,\dots,n-1,n\right\}\]
		Consequently, the cumulative distribution function (c.d.f.) of \(D_{1,n}\) is
		\begin{equation}\label{eq:3.1}
			F_{D_{1,n}}(d)=\sum_{h=1}^{d}\left(\frac{ p(1-p)^{h-1}(1-(1-p)^{n-h+1})}{1-(1-p)^{n+1}-(n+1)p(1-p)^n}\right).
		\end{equation} 
		\[
		F_{D_{1,n}}(d)=\frac{1}{1-(1-p)^{n+1}-(n+1)p(1-p)^n}\cdot\sum_{h=1}^{d} p(1-p)^{h-1}-p(1-p)^{n}.
		\]
		To simplify equation the above equation, note that
		\begin{align*}
			\sum_{h=1}^{d}p(1-p)^{h-1}-\sum_{h=1}^{d}p(1-p)^{n}
			&=p\sum_{h=0}^{d-1}(1-p)^{h}-p(1-p)^{n}\sum_{h=1}^{d}1\\
			&=p\frac{1-(1-p)^d}{p}-pd(1-p)^n.
		\end{align*}
		Hence, \[F_{D_{1,n}}(d)=\frac{1-(1-p)^d-dp(1-p)^n}{1-(1-p)^{n+1}-(n+1)p(1-p)^n}.\] It follows immediately that \[\lim_{n\to\infty}F_{D_{1,n}}(d)=1-(1-p)^d.\]
	\end{proof}

	\subsection{General Case}
	\label{sec:gc}
	
	\begin{proof}
		We once again consider the transformation \(D_{i,n}=n\Delta_{i,n}\). Using the result from Theorem \ref{thm:1} we find the c.d.f.\ of \(D_{i,n}\) is defined by \[F_{D_{i,n}}(d)=\sum_{h=1}^{d}\left[\frac{\displaystyle p^{i+1} (1-p)^{h-1}\sum_{j=i-1}^{n-h}\binom{j}{i-1}(1-p)^{j-i+1}}{\displaystyle 1-\sum_{k=0}^{i}\binom{n+1}{k}p^k(1-p)^{n+1-k}}\right],\] to evaluate its limit with respect to \(n\), we start from its numerator, that is:
		\begin{equation}\label{eq:3.2.1}
			\lim_{n\to\infty}\frac{p^{i+1}}{(1-p)^{i-1}}\sum_{h=1}^{d}\left[(1-p)^{h-1}\sum_{j=i-1}^{n-h}\binom{j}{i-1}(1-p)^{j}\right].
		\end{equation} 
		Since \(h\) is finite, the upper bound \(n-h\) of the binomial summation  will diverge as \(n\to\infty\). In addition, since \(\binom{n}{k}=0\) when \(n<k\), we can switch the lower bound of the binomial summation from \(i-1\) to simply \(0\). These allow us to safely rewrite \eqref{eq:3.2.1} as
		\begin{equation}\label{eq:3.2.2}
			\frac{p^{i+1}}{(1-p)^{i-1}}\sum_{h=1}^{d}\left[(1-p)^{h-1}\sum_{j=0}^{\infty}\binom{j}{i-1}(1-p)^{j}\right].
		\end{equation}
		To evaluate the binomial sum involved in \eqref{eq:3.2.2}, we shall utilize the method given by \autocite[p. 17]{wilf06}: let \(\langle x^n\rangle f(x)\) denote the coefficient of \(x^n\) in some power series \(f(x)\) of \(x\). Additionally, choose \(q\in\bbr\) such that \(|q|<1\) and \((1+q)(1-p)<1\), by binomial theorem we have
		\begin{equation}\label{eq:3.2.3}
			\sum_{i=0}^{\infty}\binom{j}{i-1}q^{i-1}=(1+q)^j.
		\end{equation} 
		If we multiply left hands side of the equation \eqref{eq:3.2.3} by \((1-p)^j\) and sum over \(j\geq0\), we find that \[\sum_{j=0}^{\infty}\sum_{i=0}^{\infty}\binom{j}{i-1}q^{i-1}(1-p)^j=\sum_{j=0}^{\infty}(
		(1+q)^j(1-p)^j)=\frac{1}{1-(1+q)(1-p)}.\] It follows then
		\begin{align*}
			\sum_{j=0}^{\infty}\binom{j}{i-1}(1-p)^j
			&=\langle q^{i-1}\rangle\sum_{j=0}^{\infty}\sum_{i=0}^{\infty}\binom{j}{i-1}q^{i-1}(1-p)^j\\
			&=\langle q^{i-1}\rangle\frac{1}{1-(1+q)(1-p)}\\
			&=\frac{1}{p}\langle q^{i-1}\rangle\frac{1}{1-((1-p)/p)q}\\
		\end{align*}
		Note that \[\sum_{i=1}^{\infty}\left(\frac{1-p}{p}\right)^{i-1}q^{i-1}=\frac{1}{1-((1-p)/p)q}.\] Thus,
		\begin{align*}
			\frac{1}{p}\langle q^{i-1}\rangle\frac{1}{1-((1-p)/p)q}&=\frac{1}{p}\left(\frac{1-p}{p}\right)^{i-1}\\
			&=\frac{(1-p)^{i-1}}{p^i},
		\end{align*}
		which implies \[\sum_{j=0}^{\infty}\binom{j}{i-1}(1-p)^j=\frac{(1-p)^{i-1}}{p^i}.\]
		As a result, expression \eqref{eq:3.2.2} becomes 
		\begin{equation}
			\begin{aligned}\label{eq:3.2.4}
				\frac{p^{i+1}}{(1-p)^{i-1}}&\sum_{h=1}^{d}\left[(1-p)^{h-1}\sum_{j=0}^{\infty}\binom{j}{i-1}(1-p)^{j}\right]\\&=\frac{p^{i+1}}{(1-p)^{i-1}}\sum_{h=1}^{d}\left[(1-p)^{h-1}\cdot\frac{(1-p)^{i-1}}{p^i}\right]\\
				&=p\sum_{h=1}^{d}(1-p)^{h-1}\\
				&=1-(1-p)^d.
			\end{aligned}
		\end{equation}
		Now we will examine the denominator of the c.d.f.: notice that the summation part is in fact c.d.f.\ of a binomial random variable. Through a simple transformation \(m=n+1\) we see \[\lim_{n\to\infty}\sum_{k=0}^{i}\binom{n+1}{k}p^k(1-p)^{n+1-k}=\lim_{m\to\infty}\sum_{k=0}^{i}\binom{m}{k}p^k(1-p)^{m-k}.\] Then by the De Moivre--Laplace Integral Limit Theorem \autocite[p. 33]{sinai92},
		\begin{equation}\label{eq:3.2.5}
			\lim_{m\to\infty}\sum_{k=0}^{i}\binom{m}{k}p^k(1-p)^{m-k}=\frac{1}{\sqrt{2\pi}}\int_{-\infty}^{A}e^{-z^2/2}\,dz,
		\end{equation} 
		where
		\begin{align*}
			A&=\lim_{m\to\infty}\frac{i-mp}{\sqrt{mp(1-p)}}\\&=\frac{1}{\sqrt{p(1-p)}}\lim_{m\to\infty}\frac{i-mp}{\sqrt{m}}\\
			&=\frac{-2p}{\sqrt{p(1-p)}}\lim_{m\to\infty}\sqrt{m}\\
			&=-\infty.
		\end{align*}
		Hence, equation \eqref{eq:3.2.5} yields \[\lim_{m\to\infty}\sum_{k=0}^{i}\binom{m}{k}p^k(1-p)^{m-k}=0.\] It follows then
		\begin{equation}\label{eq:3.2.6}
			\lim_{n\to\infty}1-\sum_{k=0}^{i}\binom{n+1}{k}p^k(1-p)^{n+1-k}=1.
		\end{equation}
		As a result of equations \eqref{eq:3.2.4} and \eqref{eq:3.2.6}, we conclude that \[\lim_{n\to\infty}F_{D_{i,n}}(d)=1-(1-p)^d.\] 
	\end{proof}
	
	\subsection{Alternative Proof to Theorem \ref{thm:3}}
	\label{sec:ap}
	
	We should now notice that in the case of \(n\to\infty\), our sampling process essentially becomes a infinite Bernoulli process. This allows us to investigate the asymptotic behavior of \(D_{i,n}\) from a different point of view:
	\begin{proof}
		Let \(N_k\) denote the index of the \(k\)-th survived element of \(S_n\). In this way we see that \(s'_i=s_{N_i}\). It should be noted that \(N_k\) is in fact the \(k\)-th arrival time of the process. We can then proceed to define \(M_1=N_1\) and \(M_k=N_k-N_{k-1}\) for all \(k\geq 2\) as the \(k\)-th inter-arrival time. 
		
		Suppose \(M_1=N_1=n\), then all \(n-1\) elements before \(s_{N_1}\) failed to survive. Hence, \[P(M_1=n)=(1-p)^{n-1}p.\] In other words, \(M_1\) has geometric distribution with parameter \(p\). Further, suppose that \(M_2>n\) given that \(M_1=m\). Then elements from \(s_m\) to \(s_{m+n}\) did not survive. Therefore, \[P(M_2>n\mid M_1=m)=(1-p)^n.\] This implies \(M_2\) is also of geometric distribution with parameter \(p\) and independent of \(M_1\). As a result, by using strong induction, we can establish that all \(M_k\)s are i.i.d.\ geometric random variables with parameter \(p\). 
		
		Now let's return to the scaled random spacing. Using the definition,
		\begin{align*}
			\lim_{n\to\infty}D_{i,n}&=\lim_{n\to\infty}n\Delta_{i,n}\\
			&=\lim_{n\to\infty}n\left(s'_{i+1}-s'_i\right)\\
			&=\lim_{n\to\infty}n\left(s_{N_{i+1}}-s_{N_i}\right)\\
			&=\lim_{n\to\infty}n\left(\frac{N_{i+1}-1}{n}-\frac{N_{i}-1}{n}\right)\\
			&=M_{i+1}.
		\end{align*} 
		It follows immediately that \(\{D_{i,n}\}\) converges in distribution to a geometric random variable with parameter \(p\), as demonstrated in previous section.
	\end{proof}
	
	\section*{Acknowledgements}
	We would like to thank Professor Jayadev Athreya and Professor Francesco Cellarosi, the faculty advisors of the research. We would like to thank Professor Trevor Park for probability theory related discussion. We would like to thank Amita Malik and Elizabeth C. Merriman, the graduate mentors of the research. We would also like to thank Suhan Liu, former member of the research team.
	\newpage
	\printbibliography
	
\end{document}